\documentclass[12pt,reqno]{amsart}
\usepackage{amssymb}

\newtheorem{theorem}{Theorem}[section]

\theoremstyle{definition}

\numberwithin{equation}{section}

\begin{document}

\baselineskip=15.5pt

\title[Quot schemes and Ricci semipositivity]{Quot schemes and Ricci semipositivity}

\author[I. Biswas]{Indranil Biswas}

\address{School of Mathematics, Tata Institute of Fundamental
Research, Homi Bhabha Road, Mumbai 400005, India}

\email{indranil@math.tifr.res.in}

\author[H. Seshadri]{Harish Seshadri}

\address{Indian Institute of Science, Department of Mathematics,
Bangalore 560003, India}

\email{harish@math.iisc.ernet.in}

\subjclass[2000]{14H60, 14H81, 32Q10}

\date{}

\begin{abstract}
Let $X$ be a compact connected Riemann surface of genus at least two, and let
${\mathcal Q}_X(r,d)$ be the quot scheme that parametrizes all the
torsion coherent quotients of ${\mathcal O}^{\oplus r}_X$ of degree $d$.
This ${\mathcal Q}_X(r,d)$ is also a moduli space of vortices on $X$. Its geometric
properties have been extensively studied. Here we prove that the anticanonical line bundle
of ${\mathcal Q}_X(r,d)$ is not nef. Equivalently, ${\mathcal Q}_X(r,d)$ does not admit any
K\"ahler metric whose Ricci curvature is semipositive.
\medskip

\textbf{R\'esum\'e}. \textit{Sch\'ema quot et semi-positivit\'e de Ricci}\\

Soit $X$ une surface de Riemann compacte et connexe de genre au moins deux, et soit 
${\mathcal Q}_X(r,d)$ le sch\'ema quot qui param\'etrise tous les quotients torsion 
coh\'erents de ${\mathcal O}^{\oplus r}_X$ de degr\'e $d$. L'espace ${\mathcal Q}_X(r,d)$ est 
aussi un espace de modules de vortex sur $X$. Nous d\'emontrons que le fibr\'e anticanonique 
de $X$ n'a pas la propri\'et\'e nef. De façon \'equivalente, ${\mathcal Q}_X(r,d)$ n'admet 
aucune m\'etrique kahl\'eienne dont la courbure de Ricci est semi-positive.
\end{abstract}

\maketitle

\section{Introduction}

Take a compact connected Riemann surface $X$. The genus of $X$, which will be denoted by $g$, is 
assumed to be at least two. We will not distinguish between the holomorphic vector bundles on $X$ 
and the torsion-free coherent analytic sheaves on $X$. For a positive integer $r$, let ${\mathcal 
O}^{\oplus r}_X$ be the trivial holomorphic vector bundle on $X$ of rank $r$. Fixing a positive 
integer $d$, let
\begin{equation}\label{e1}
{\mathcal Q}\, :=\, {\mathcal Q}_X(r,d)
\end{equation}
be the quot scheme that parametrizes all (torsion) coherent quotients of ${\mathcal
O}^{\oplus r}_X$ of rank zero and degree $d$ \cite{Gr}. Equivalently, ${\mathcal Q}$ parametrizes all
coherent subsheaves of ${\mathcal O}^{\oplus r}_X$ of rank $r$ and degree $-d$, because these
are precisely the kernels of coherent quotients of ${\mathcal
O}^{\oplus r}_X$ of rank zero and degree $d$. This $\mathcal Q$ is a connected smooth
complex projective variety of dimension $rd$. See \cite{BGL}, \cite{Bi}, \cite{BDW} for
properties of $\mathcal Q$. It should be mentioned that $\mathcal Q$ is also
a moduli space of vortices on $X$, and it has been extensively studied from this point of
view of mathematical physics; see \cite{Ba}, \cite{BiRo}, \cite{BR} and references therein.

B\"okstedt and Rom\~{a}o proved some interesting differential geometric properties of
$\mathcal Q$ (see \cite{BR}). In \cite{BS0} and \cite{BS} we proved that $\mathcal Q$
does not admit K\"ahler metrics with semipositive or seminegative holomorphic
bisectional curvature. In this note, we continue the study the question of existence of
metrics on $\mathcal Q$ whose curvature has a sign. Our aim here is to
prove the following:

\begin{theorem}\label{thm1}
The quot scheme ${\mathcal Q}$ in \eqref{e1} does not admit any K\"ahler metric such
that the anticanonical line bundle $K^{-1}_{\mathcal Q}$ is hermitian semipositive.
\end{theorem}

Since semipositive holomorphic bisectional curvature implies semipositive Ricci
curvature for a K\"ahler metric, Theorem \ref{thm1} generalizes the main result of
\cite{BS}.

Recall that a holomorphic line bundle $L$ on a compact complex manifold $M$ is said to
be \textit{hermitian semipositive} if $L$ admits a smooth hermitian structure such
that the corresponding hermitian connection has the property that its curvature form
is semipositive. The anticanonical line bundle on $M$ will be denoted by $K^{-1}_M$.
Note that if $M$ admits a K\"ahler metric such that the corresponding Ricci curvature
is semipositive, then $K^{-1}_M$ is hermitian semipositive. Indeed, in that case the
hermitian connection on $K^{-1}_M$ for the hermitian structure induced by such a
K\"ahler metric has semipositive curvature. The converse statement, that hermitian
semipositivity of $K_M^{-1}$ implies the existence of K\"ahler metrics with
semipositive Ricci curvature, is also true by Yau's solution of the Calabi's conjecture
\cite{aub1}, \cite{aub2}, \cite{yau}.

The proof of Theorem \ref{thm1} is based on a recent work of Demailly, Campana and Peternell on 
the classification of compact K\"ahler manifolds $M$ with semipositive $K_M^{-1}$ \cite{De}, 
\cite{CDP}. This classification implies that if $K_M^{-1}$ is semipositive, then there is a 
nontrivial abelian ideal in the Lie algebra of holomorphic vector fields on $M$, provided 
$b_1(M)\, >\, 0$. On the other hand, for $M\,=\, {\mathcal Q}$, this Lie algebra is isomorphic to 
$\mathfrak {sl}(r, {\mathbb C})$, which does not have any nontrivial abelian ideal.

\section{Proof of Theorem \ref{thm1}}

\subsection{Semipositive Ricci curvature}

Let $J^d(X)\,=\, \text{Pic}^d(X)$ be the connected component of the Picard group of $X$ that
parametrizes the isomorphism classes of holomorphic line bundles on $X$ of degree $d$. Let
$S^d(X)$ denote the space of all effective divisors on $X$ of degree $d$, so
$S^d(X)\,=\, X^d/P_d$ is the symmetric product with $P_d$ being the group of permutations of
$\{1\, ,\cdots\, ,d\}$. Let
\begin{equation}\label{p}
p\, :\, S^d(X)\,\longrightarrow\, \text{Pic}^d(X)
\end{equation}
be the natural morphism that sends a divisor on $X$ to the holomorphic line bundle
on $X$ defined by it.

Take any coherent subsheaf $F\,\subset\, {\mathcal O}^{\oplus r}_X$ of
rank $r$ and degree $-d$. Let
$$
s_F\, :\, {\mathcal O}^{\oplus r}_X\,=\, ({\mathcal O}^{\oplus r}_X)^*\,
\longrightarrow\, F^*
$$
be the dual of the inclusion of $F$ in ${\mathcal O}^{\oplus r}_X$. Its exterior
product
$$
\bigwedge\nolimits^r s_F\, :\, {\mathcal O}_X\,=\, \bigwedge\nolimits^r{\mathcal O}^{\oplus r}_X\,
\longrightarrow\, \bigwedge\nolimits^r F^*
$$
is a holomorphic section of the holomorphic line bundle $\bigwedge^r F^*$ of degree $d$.
Therefore, the divisor $\text{div}(\bigwedge^r s_F)$ is an element of $S^d(X)$. Consequently,
we have a morphism
\begin{equation}\label{e2}
\varphi\, :\, {\mathcal Q}\,\longrightarrow\, S^d(X)\, ,\ \ F\, \longmapsto\,
\text{div}(\bigwedge\nolimits^r s_F)\, ,
\end{equation}
where $\mathcal Q$ is defined in \eqref{e1}. We note that when $r\,=\,1$, then
$\varphi$ is an isomorphism.

Assume that ${\mathcal Q}$ admits a K\"ahler metric $\omega$ such that $K^{-1}_{\mathcal
Q}$ is hermitian semipositive. Then there is a connected finite \'etale Galois covering
\begin{equation}\label{f}
f\, :\, \widetilde{\mathcal Q}\, \longrightarrow\, \mathcal Q
\end{equation}
such that $(\widetilde{\mathcal Q},\, f^*\omega)$ is holomorphically isometric to a product
\begin{equation}\label{e3}
\gamma\, :\, \widetilde{\mathcal Q}\, \longrightarrow\, A\times C\times H\times F\, ,
\end{equation}
where
\begin{itemize}
\item{} $A$ is an abelian variety,

\item $C$ is a simply connected Calabi--Yau manifold (holonomy is $\text{SU}(c)$, where $c\,
=\, \dim C$),

\item $H$ is a simply connected hyper-K\"ahler manifold (holonomy is $\text{Sp}(h/2)$, where
$h\, =\, \dim H$), and

\item{} $F$ is a rationally connected smooth projective variety such that $K^{-1}_F$ is
hermitian semipositive.
\end{itemize}
(See \cite[Theorem~3.1]{De}.) Henceforth, we will identify $\widetilde{\mathcal Q}$ with $A\times 
C\times H\times F$ using $\gamma$ in \eqref{e3}. We note that $F$ is simply connected because it 
is rationally connected \cite[p.~545, Theorem 3.5]{Ca}, \cite[p.~362, Proposition 2.3]{Ko}.

\subsection{A lower bound of $d$}

We know that $b_1({\mathcal Q})\,=\, 2g$, and the induced homomorphism
$$
(p\circ\varphi)_*\, :\, H_1({\mathcal Q},\, {\mathbb Q})\, \longrightarrow\,
H_1(\text{Pic}^d(X),\, {\mathbb Q})\, ,
$$
where $p$ and $\varphi$ are constructed in \eqref{p} and \eqref{e2} respectively, is
an isomorphism \cite{Bi}, \cite[p.~649, Remark]{BGL}.
Since $f$ in \eqref{f} is a finite \'etale covering, the induced homomorphism
$$
f_*\, :\, H_1(\widetilde{\mathcal Q},\, {\mathbb Q})\, \longrightarrow\,
H_1({\mathcal Q},\, {\mathbb Q})
$$
is surjective. Therefore, the homomorphism
\begin{equation}\label{e5}
(p\circ\varphi\circ f)_*\, :\, H_1(\widetilde{\mathcal Q},\, {\mathbb Q})\, \longrightarrow\,
H_1(\text{Pic}^d(X),\, {\mathbb Q})
\end{equation}
is surjective.

There is no nonconstant holomorphic map from a compact simply connected K\"ahler
manifold to an abelian variety. In particular, there are no nonconstant holomorphic maps from
$C$, $H$ and $F$ in \eqref{e3} to $\text{Pic}^d(X)$. Hence the map
$p\circ\varphi\circ f$ factors through a map
$$
\beta\, :\, A\, \longrightarrow\, \text{Pic}^d(X)\, .
$$
In other words, there is a commutative diagram
\begin{equation}\label{cd}
\begin{matrix}
\widetilde{\mathcal Q} = A\times C\times H\times F &\
\stackrel{p\circ\varphi\circ f}{\longrightarrow} & \text{Pic}^d(X)\\
q\Big\downarrow~{\ } && ~\,~\Vert{\rm Id}\\
A &\ \stackrel{\beta}{\longrightarrow} & \text{Pic}^d(X)
\end{matrix}
\end{equation}
where $q$ is the projection of $A\times C\times H\times F$ to the first factor. Since
$H_1(A\times C\times H\times F,\, {\mathbb Z})\,=\, H_1(A,\, {\mathbb Z})$ (as $C$, $H$ and
$F$ are simply connected), and $(p\circ\varphi\circ f)_*$ in \eqref{e5} is surjective, it
follows that the homomorphism
$$
\beta_*\, :\, H_1(A,\, {\mathbb Q}) \, \longrightarrow\,H_1(\text{Pic}^d(X),\, {\mathbb Q})
$$
induced by $\beta$ is surjective. This immediately implies that the map $\beta$ is
surjective. Since $\beta$ is surjective, from the commutativity of \eqref{cd}
we know that the map $p$ is surjective. This implies that
\begin{equation}\label{ge}
d\, =\, \dim S^d(X) \, \geq\, \dim \text{Pic}^d(X)\,=\, g\, \geq\, 2\, .
\end{equation}

\subsection{Albanese for $\widetilde{\mathcal Q}$}

The homomorphism of fundamental groups
$$
\varphi_*\, :\, \pi_1({\mathcal Q})\, \longrightarrow\, \pi_1(S^d(X))
$$
induced by $\varphi$ in \eqref{e2} is an isomorphism \cite[Proposition 4.1]{BDHW}. Since
$d\, \geq\, 2$ (see \eqref{ge}), the homomorphism of fundamental groups
$$
p_*\, :\, \pi_1(S^d(X))\, \longrightarrow\,\pi_1(\text{Pic}^d(X))
$$
induced by $p$ in \eqref{p} is an isomorphism. Indeed, $\pi_1(S^d(X))$ is the
abelianization $$\pi_1(X)/[\pi_1(X),\, \pi_1(X)]\,=\, H_1(X, \, {\mathbb Z})$$ of
$\pi_1(X)$ \cite{DT}. Combining these we conclude that the homomorphism of fundamental groups
\begin{equation}\label{pi1}
(p\circ\varphi)_*\, :\, \pi_1({\mathcal Q})\, \longrightarrow\, \pi_1(\text{Pic}^d(X))
\end{equation}
induced by $p\circ\varphi$ is an isomorphism.

Since the homomorphism in \eqref{pi1} is an isomorphism, the covering $f$ in \eqref{f} is induced
by a covering of $\text{Pic}^d(X)$. In other words, there is a finite \'etale Galois covering
\begin{equation}\label{mu}
\mu\, :\, J\, \longrightarrow\, \text{Pic}^d(X)
\end{equation}
and a morphism $\lambda\, :\, \widetilde{\mathcal Q}\,\longrightarrow\, J$ such that the following
diagram is commutative:
\begin{equation}\label{cc}
\begin{matrix}
\widetilde{\mathcal Q} &\ \stackrel{f}{\longrightarrow} & {\mathcal Q}\\
{\ }~\Big\downarrow\lambda && {\ \ \ \ \ \ }~~~\Big\downarrow p\circ\varphi\\
J &\ \stackrel{\mu}{\longrightarrow} & \text{Pic}^d(X)
\end{matrix}
\end{equation}
where $f$ is the covering map in \eqref{f}. The projection $q$ in \eqref{cd} is clearly the 
Albanese morphism for $\widetilde{\mathcal Q}$, because $C$, $H$ and $F$ are all simply 
connected. On the other hand, $p\circ\varphi$ is the Albanese morphism for $\mathcal Q$ 
\cite[Corollary 2.2]{BS}. Therefore, its pullback, namely, $\lambda$, is the Albanese morphism 
for $\widetilde{\mathcal Q}$. Consequently, we have $A\,=\, J$ with $\lambda$ coinciding with the 
projection $q$ in \eqref{cd}. Henceforth, we will identify $A$ and $q$ with $J$ and $\lambda$ 
respectively.

\subsection{Vector fields}

The differential $df$ of $f$ identifies $T{\widetilde{\mathcal Q}}$ with $f^*T{\mathcal Q}$,
because $f$ is \'etale.
Using the trace homomorphism $t\, :\, f_*{\mathcal O}_{\widetilde{\mathcal Q}}\, \longrightarrow\,
{\mathcal O}_{\mathcal Q}$, we have
$$
f_*T{\widetilde{\mathcal Q}}\,=\, f_*f^*T{{\mathcal Q}}\, \stackrel{p_f}{\longrightarrow}\,
(f_*{\mathcal O}_{\widetilde{\mathcal Q}})\otimes T{{\mathcal Q}}\,
\stackrel{t}{\longrightarrow}\,{\mathcal O}_{\mathcal Q}\otimes T{{\mathcal Q}}\,=\,
T{{\mathcal Q}}\, ,
$$
where $p_f$ is given by the projection formula. This produces a homomorphism
\begin{equation}\label{Phi}
\Phi\, :\, H^0(\widetilde{\mathcal Q},\, T{\widetilde{\mathcal Q}})\,=\, H^0({\mathcal Q},\,
f_*T{\widetilde{\mathcal Q}})\,\longrightarrow\, H^0({\mathcal Q},\, T{{\mathcal Q}})
\end{equation}
(the equality $H^0(\widetilde{\mathcal Q},\, T{\widetilde{\mathcal Q}})\,=\, H^0({\mathcal Q},\,
f_*T{\widetilde{\mathcal Q}})$ follows from the fact that $f$ is a finite morphism).
This homomorphism $\Phi$ is surjective. Indeed, as $f^*T{\mathcal Q}\,=\,T{\widetilde{\mathcal
Q}}$, any section of $T{{\mathcal Q}}$ pulls back to a section of $T{\widetilde{\mathcal Q}}$.

Since $\widetilde{\mathcal Q}\, =\, A\times C\times H\times F$, we have
\begin{equation}\label{Phi2}
H^0(\widetilde{\mathcal Q},\, T{\widetilde{\mathcal Q}})\,=\,
H^0(A,\, TA)\oplus H^0(C,\, TC)\oplus H^0(H,\, TH)
\oplus H^0(F,\, TF)\, .
\end{equation}
Note that $H^0(\widetilde{\mathcal Q},\, T{\widetilde{\mathcal Q}})$ is a Lie algebra
under the operation of Lie bracket of vector fields, and the subspace
$$
H^0(A,\, TA)\, \subset\, H^0(\widetilde{\mathcal Q},\, T{\widetilde{\mathcal Q}})
$$
(see \eqref{Phi2}) is an ideal in this Lie algebra. Since $A\,=\, J$ is a covering of
$\text{Pic}^d(X)$, we have
\begin{equation}\label{dg}
\dim H^0(A,\, TA)\,=\, \dim \text{Pic}^d(X)\,=\, g\, >\, 1\, .
\end{equation}

Since $H^0(A,\, TA)$ is an ideal in $H^0(\widetilde{\mathcal Q},\, T{\widetilde{\mathcal
Q}})$, it follows immediate that
$$
\Phi(H^0(A,\, TA))\, \subset\, \Phi(H^0(\widetilde{\mathcal Q},\, T{\widetilde{\mathcal Q}}))\,=\,
H^0({\mathcal Q},\, T{\mathcal Q})
$$
is an ideal, where $\Phi$ is constructed in \eqref{Phi}. Note that $H^0(A,\, TA)$ is an abelian
Lie algebra, so the Lie algebra $\Phi(H^0(A,\, TA))$ is also abelian.

Since $\mu\, :\, J\,=\, A\, \longrightarrow\, \text{Pic}^d(X)$ in \eqref{mu} is a covering map
between abelian varieties, the trace map $H^0(A,\, TA) \, \longrightarrow\,
H^0(\text{Pic}^d(X),\, T\text{Pic}^d(X))$ is an isomorphism. In view of this, from the
commutativity of the diagram in \eqref{cc} it follows that the restriction
$$
\Phi\vert_{H^0(A,\, TA)}\, :\, H^0(A,\, TA)\, \longrightarrow\, H^0({\mathcal Q},\, T{\mathcal Q})
$$
is injective (see \eqref{Phi2} and \eqref{Phi}). But $H^0({\mathcal Q},\, T{\mathcal Q})\,=\, 
\mathfrak{sl}(r, {\mathbb C})$ \cite[p.~1446, Theorem 1.1]{BDH}. Hence the Lie algebra 
$H^0({\mathcal Q},\, T{\mathcal Q})$ does not contain any nonzero abelian ideal. This is in 
contradiction with the earlier result that $\Phi(H^0(A,\, TA))$ is a nonzero abelian ideal in 
$H^0({\mathcal Q},\, T{\mathcal Q})$ of dimension $g$ (see \eqref{dg}). This completes the proof 
of Theorem \ref{thm1}.


\end{document}